\documentclass[11pt]{article}
\usepackage{amsmath, amsthm, amssymb}    
\usepackage{fullpage}
\usepackage{graphicx}                    
\usepackage{hyperref}                    

\usepackage{multicol}

\title{On Mathematical Ways of Knowing}
\author{Gizem Karaali\\
Dept.~of Mathematics, Pomona College\\
610 North College Avenue, Claremont, CA 91711, USA\\
{\tt gizem.karaali@pomona.edu}
}
\date{}				

\begin{document}

\maketitle

\thispagestyle{empty}

%
%



%


Humanistic mathematics is a perspective on mathematics that emphasizes the ways our species creates, interacts with, and lives through it. I summarized this idea elsewhere (see \cite{Karaali15}) by asserting that mathematics is the way our species makes sense of this world and that it is inherent in our thinking machinery; our mathematics in turn is dependent on the way we view our universe and ourselves. Lakoff and N\'{u}\~{n}ez \cite{Lakoff} argue carefully and eloquently for a mathematics inherently based on human cognition. 

Cognition is ``the mental action or process of acquiring knowledge and understanding through thought, experience, and the senses" (see \href{https://en.wikipedia.org/wiki/Cognition}{Wikipedia}). In this note I attempt to engage with the construct of mathematical cognition through the lens of humanistic mathematics.

\vspace{-10pt} \section*{three questions}

Cognition is essentially about mental processes involving knowledge, knowing, and understanding; mathematical cognition therefore raises questions about mathematical knowledge, knowing mathematics, and understanding mathematics. Thus, I first intend to explore broadly three related questions: 

\begin{description}
\item[1.] What does it mean to know something mathematical? 
\item[2.] How do we come to know a mathematical truth?
\item[3.] What does it mean to understand mathematics? 
\end{description}

In what follows, I will not pretend to offer a comprehensive treatment of any of these questions. But in the very least I intend to open up all three questions in productive ways, so that all readers intrigued by these questions will find in the following assertions worth agreeing with and arguing against. 

\vspace{-10pt} \section*{question 1}

The first question is a natural extension of traditional epistemological investigations into mathematics. Philosophers have tinkered with the knowledge question for centuries, or rather, millennia, and mathematical knowledge has often been a part of the equation. Knowledge as justified true belief, a core tenet of epistemology since the Enlightenment, is where I want to start this note.\footnote{  At least since the 1960s with the publication of Edmund Gettier's ``Is Justified True Belief Knowledge?" \cite{Gettier}, this conventional approach to knowledge has seen many rebuttals and rephrasals, but for this note, I will mostly ignore this recent body of work. My interest is most in line with the idea of mathematical knowledge as being justified true belief.}

If (mathematical) knowledge is justified true belief (in mathematical statements), we have multiple avenues to the first question. Or alternatively we have two related questions to attend to:

\begin{description}
\item[1A.] What does it mean that a mathematical assertion is true?
\item[1B.] What does it mean that a belief in a mathematical assertion is justified?
\end{description}

1A is perhaps on the natural playground of mathematicians. Mathematicians seem to be concerned quite single-mindedly and profoundly in the truth of their assertions. One can even suggest that mathematics is nothing if not true. That is, doing math is making true mathematical assertions.\footnote{  Truth is surely not the only target of mathematicians. It is not even the driving force, according to William Byers, who writes: ``Classifying ideas as true or false is just not the best way of thinking about them. Ideas may be fecund; they may be deep; they may be subtle; they may be trivial. These are the kinds of attributes we should ascribe to ideas." \cite{Byers}. As a working mathematician, I agree with Byers, but this does not mean that I don't also believe that truth is a prerequisite. Even when mathematicians work with tentative and even patently false assertions, they have a broader truth in mind, and are not done until eventually they can reach that truth.}  In some sense, therefore, I think that the truth of a mathematical assertion means that it is a part of mathematics, this edifice we human mathematicians are building together. Philosophers have tried to clarify what mathematicians might mean when they say that a mathematical statement is true. Sitting between mathematics and philosophy, the logician Alfred Tarski \cite{Tarski3356} proposed a definition of just what truth might mean in ``the deductive sciences", which presumably include mathematics. Once again many have commented on Tarski's definition of truth. I will not go into that here but there is indeed much more that can be said along these lines if one is concerned about question 1A.

If we want to address question 1B about justification, we can, like some, invoke idealized conceptions of mathematical justification involving formal systems and proof theory. 
But most mathematicians agree that belief in the truth of a mathematical statement is justified once there is a proof of the statement that experts can agree upon. This is quite in tune with Reuben Hersh's various definitions of a ``proof" in \cite{Hersh}, most notably ``The `proof' is a procedure, an argument, a series of claims, that every qualified expert understands and accepts." Though some philosophers reading Hersh might disagree (see for instance \cite{Pollard}), it is indeed the case that when mathematicians claim that a statement is true, they mean that there is some consensus among the relevant experts that the statement is true. And an argument might have been a proof at a given time and place and afterwards, with contemporary expertise changing sides, it might become invalid. Similarly, proposed proofs do not become proofs until verified and validated by experts. Indeed, one might argue that ``it is the provision of [\ldots] evidence, not the endorsement of experts, that makes [a display of symbols, words, diagrams and such] a proof'' \cite{Pollard}. However, nothing needs to change in the display for an argument to remain an alleged proof until experts deem it is valid, and then and only then is the rest of the mathematical community comfortable in feeling justified to believe that the mathematical statement in consideration has finally been proved. What counts as persuasive evidence is almost always context-dependent. In the case of law this is obvious; even Wikipedia knows that there are variations on what counts as proof, what counts as evidence.\footnote{See \url{https://en.wikipedia.org/wiki/Burden_of_proof_(law)} for a selection of legal standards of evidence and proof. }  Why do we expect mathematics to be different? 

If we see mathematics as something done by humans, the formalist, proof theory-based understandings of proof and mathematics remain idealized approximations at best. It is the human (or sometimes, and begrudgingly, human-assisted) verification that mathematicians look for in a proof.\footnote{  Do we know if the Four-Color Theorem is true? Yes, we do. Or at least most mathematicians would concede that the human-assisted computer proof (or alternatively, the computer-assisted human proof) is enough for us to agree that it is true. There are still those who want more human proofs of the result, but the truth of the statement does not need further justification. }  And this is definitely context-based, both in terms of space-time coordinates and the cultural makeup of the audience. As Israel Kleiner quotes in \cite{Kleiner}, G.F. Simmons wrote ``Mathematical rigor is like clothing: in its style it ought to suit the occasion, and it diminishes comfort and restricts freedom of movement if it is either too loose or too tight.'' This quote captures well how our understanding of just what should count as proof is dependent on the fashions of our times. This aligns with Harel's perspective \cite{Harel}: ``Mathematics is a human endeavor, not a predetermined reality. As such, it is the community of the creators of mathematics who makes decisions about the inclusion of new discoveries in the existing edifice of mathematics.'' Among the decisions left to the human creators of mathematics are the truth of a mathematical statement and the validity of its proof. 

\vspace{-10pt} \section*{question 2} 

The standard modern answer to question 2 is ``by a rigorous proof". Let us leave aside the historicity and cultural dependency of this response now, and its vagueness (what is proof? what is rigor?). I already wrote a bit about all that above. The reader who is not yet convinced may read \cite{Kleiner} for more on rigor and proof.\footnote{A tangentially related read which might add some insight to this conversation is \cite{Tao} about the role of rigor and proof in mathematics.} But I want to emphasize here the possible distinction between the doer of the proof and the believer who believes with justification that the proof is valid and that the related statement is true. Does the believer need to understand the proof in order to know the related statement is true? How similar and how different is this from the calculus student saying that they learned calculus because they passed the final exam? Let us narrow things down a bit more. Should the successful calculus student be able to state the Fundamental Theorem of Calculus? Should they be able to prove it? Should they be able to replicate the argument in their textbook or the one their instructor put on the board? Should they be able to provide a convincing argument for its truth? Alternatively, should they be able to use it in a problem that requires the result? When do we assume a student has learned, or knows the Fundamental Theorem of Calculus?

Mathematics education researcher Guershon Harel thinks that these kinds of pedagogical questions are not independent from the philosopher's concern about mathematical knowledge. In particular he proposes a definition of mathematics which originates from his pedagogical research that might help us with our endeavor here:
\begin{quote}
Mathematics consists of two complementary subsets:\\
$\bullet$ The first subset is a collection, or structure, of structures consisting of particular axioms, definitions, theorems, proofs, problems, and solutions. This subset consists of all the institutionalized ways of understanding in mathematics throughout history. It is denoted by WoU.\\
$\bullet$ The second subset consists of all the ways of thinking, which are characteristics of the mental acts whose products comprise the first set. It is denoted by WoT. \cite{Harel}.
\end{quote}

Here Harel uses ``ways of understanding'' and ``ways of thinking'' as technical terms. According to him ``a way of understanding is a particular cognitive product of a mental act carried out by an individual'', while ``a way of thinking, on the other hand, is a cognitive characteristic of a mental act''. Here is how he uses them in context:
\begin{quote}
``Any statement a teacher (or a classmate) utters or puts on the board will be translated by each individual student into a way of understanding that depends on her or his experience and background.'' 

``The range of ways of understanding a fraction makes the area of fractions a powerful elementary mathematics topic?one that can offer young students a concrete context to construct desirable?indeed, crucial?ways of thinking, such as: mathematical concepts can be understood in different ways, mathematical concepts should be understood in different ways, and it is advantageous to change ways of understanding of a mathematical concept in the process of solving problems.''
\end{quote}
Consider the mathematician who reads the statement of the Four-Color Theorem and discussions about its proof and as a result is convinced of the veracity of the statement and the sufficiency of the proof. This mathematician, in my opinion, is not that different from the good student who learned of the Fundamental Theorem of Calculus from their instructor, can use it in various scenarios, can even perhaps outline a convincing argument about why it might be a reasonable thing to assume. In each case I'd say the person knows the concept and construct in question. They believe the truth of a true statement and are justified in doing so. They put their trust in a relatively trustworthy source of expertise. But does the student really know (that is, understand) the Fundamental Theorem of Calculus? Does the mathematician really know (that is, understand) the Four-Color Theorem? This naturally brings us to question 3: What does it mean to understand something in mathematics?

\vspace{-10pt} \section*{question 3}

Harel's ways of thinking and ways of understanding are related quite visibly to question 3. For example, Harel offers a handful of ways of understanding the concept of fractions: 
\begin{enumerate}
\item[a.] the part-whole interpretation: $m / n$ (where $m$ and $n$ are positive integers) means ``$m$ out of $n$ objects.'' 
\item[b.] $m/n$ means ``the sum $\tfrac{1}{n} + \cdots + \tfrac{1}{n}$, $m$ times''
\item[c.] ``the quantity that results from $m$ units being divided into $n$ equal parts'' 
\item[d.] ``the measure of a segment $m$ inches long in terms of a ruler whose unit is $n$ inches'' 
\item[e.] ``the solution to the equation $nx = m$'' 
\item[f.] ``the ratio $m : n$; namely, $m$ objects for each $n$ objects.'' 
\end{enumerate}

Similarly, we can develop a list of ways of understanding for the derivative, in terms of the limit definition; in terms of slopes of tangent lines; in terms of linear approximations; and so on. And it seems reasonable to assume that when we say that a student understands fractions, we mean that they have mastered an indeterminate (but definitely nonzero) number of ways of understanding the concept. The fluency with which they can move from one interpretation to the other can help us if we want to further qualify how much they understand. 

Understanding seems to presuppose knowing but is there anything more to it? More specifically, understanding a piece of mathematics does presuppose knowing that piece of mathematics; what we want to know is if it involves anything more. 

This is the perfect context for me to bring up my favorite quote from one of my mathematical heroes. John von Neumann was a polymath, a genius mathematician who was instrumental in the development of quantum mechanics, game theory, functional analysis, operator algebras, and computer science, a great mathematical mind. This distinguished mathematician is known to have said to a young scholar asking for advice: ``Young man, in mathematics you do not understand things. You just get used to them."\footnote{ For an interesting discussion of this quote see \url{https://math.stackexchange.com/questions/11267/what-are-some-interpretations-of-von-neumanns-quote}.}

Grant that this is a sharp quote. It hits you hard and shakes you up, especially if you have at least a passing knowledge of the extent of von Neumann's own mathematical contributions. But can you take it seriously? Can you get anything out of it in the context of mathematical knowledge and mathematical ways of knowing?

My personal take on this quote is twofold. 

One is that of the optimistic student of mathematics. Even if I feel like I am not understanding something, there is some benefit to pushing forward, if only a bit more. Doubtless the great mathematician is right; sometimes you simply have to move on, after accepting the fact as a fact and see where it leads you. Stubborn patience. Dogged perseverance. 

However perhaps von Neumann did not really mean to recommend moving forward without understanding. Perhaps he was saying something else, that what you call understanding is not something subtle or sublime. In fact, when we are learning a new concept, a new theory, don't we start by making mental patterns, charting new pathways in our mind, formatting our minds so certain types of programs run well or smoothly enough? How is this different cognitively from getting used to brushing our teeth before going to bed, splashing our faces before leaving the bathroom, eating with the fork in our right hand? Habit forming is done by doing something over and over again; aren't mathematical ways of thinking and ways of understanding reinforced by repeated practice as well? And is there a genuinely different, a genuinely distinct sense of ``understanding'' that goes beyond ``getting used to thinking of the concept in question in a particularly productive manner''? 

Imagine a student who learns to think of a complex number first as an ordered pair, then as a point on the complex plane, and then as a linear transformation on the complex plane. When can we say the student knows complex numbers? When do we say the student understands them? I agree with Emily Grosholz who argues, using complex numbers as a concrete case study, that ``the best way to teach students mathematics is through a repertoire of modes of representation, which is also the best way to make mathematical discoveries'' \cite{Grosholz}. But this also makes things a lot more complicated. If there are multiple ways of understanding, a la Harel, that need to go into understanding a construct, when do we really know the construct? When do we understand it?

Bloom's taxonomy is just one of many ranking frameworks education researchers use to delineate cognitive tasks and their demands on a learner. In the original taxonomy of Bloom and colleagues, knowledge is the very lowest level of cognition needed and includes knowledge of specifics, terminology, specific facts, ways and means of dealing with specifics, conventions, trends and sequences, classifications and categories, criteria, methodology, the universals and abstractions in a field, principles and generalizations, theories and structures \cite{Bloom}. Understanding is related to the second level, comprehension. (In the revised version \cite{Anderson}, to understand is once again the second level, and there are several other layers of cognition ranked higher.) Mathematics education researchers prefer other frameworks to evaluate the cognitive demand of mathematical tasks; see \cite{Smith} for a commonly used schema distinguishing between memorization, procedures without or with connections, and ``doing mathematics''. It is quite interesting that learning mathematics involves doing mathematics as a subcategory!

Putting taxonomies aside, knowing and understanding a concept might not be that different from one another after all. There seems to be a psychological difference for sure; understanding always occurs with knowing but sometimes we might feel we know something but do not ``really understand''. But perhaps von Neumann was not just being cheeky. Perhaps there is really nothing more than getting used to knowing something. Perhaps all we need is a brain formatted in the right way to accept our knowledge as truth and naturally so.\footnote{  The formatting of the brain might sound strange but is not that far from the point of education. The point of education is to shape students' minds. What is that if not brain formatting? Indoctrination might also fit in this view. Some might take this to discussions of how mathematics education promotes dogmatic beliefs. Though I am interested in such inquiry, I will not pursue it here. }

\vspace{-10pt} \section*{mathematical ways of knowing}

In the remainder of this note, I want to delineate a construct I will call ``mathematical ways of knowing''. This is in some sense related to what Harel calls ``ways of thinking''. But I believe it is not exactly the same. 

I start by the axiomatic definition that mathematics is one of the main systematic bodies of knowledge that formulates and occasionally aims to address questions about human perception of the world and the human endeavors to understand it. The concepts and constructs of number, shape, form, time, change, and chance are fundamental to our understanding of our world as well as ourselves. These concepts and constructs are mathematical in nature, or at least they are naturally amenable to mathematical approaches. 

Within this setting, I mean by a mathematical way of knowing a way of formulating and addressing a question or a set of questions about our world and ourselves that allows for mathematical inquiry. It might be interesting to try and put these mathematical ways of knowing in contrast to or in conversation with a handful of other ways of knowing: scientific / empirical, faith-based, philosophical. I do no such thing in this note however. Here I merely put down some ideas as placeholders, as tentative yet suggestive notes to a self that might or might not be able to come back to revisit them in a possible future. 

\vspace{-10pt} \section*{mathematical ways of knowing: rationalism and imagination}

The first two mathematical ways of knowing I would like to consider are rationalism and imagination. 

Rationalism in mathematics is the fundamental assumption that we ought to reach our mathematical truths through reason. Experiment and experience might provide hints towards a truth, but they are never enough to convince us in the final count. Even if we ``know something in our guts'', we are not convinced we have mathematical certainty until we can reason our way to that something. 

And reason and rational thinking are captured effectively in axiomatic thinking. The central tenets of axiomatic thinking are captured by Tarski in the following:
\begin{quote}
When we set out to construct a given discipline, we distinguish, first of all, a certain small group of expressions of this discipline that seem to us to be immediately understandable; the expressions in this group we call primitive terms or undefined terms, and we employ them without explaining their meanings. At the same time we adopt the principle: not to employ any of the other expressions of the discipline under consideration, unless its meaning has first been determined with the help of primitive terms and of such expressions of the discipline whose meanings have been explained previously. \cite{Tarski46}
\end{quote}

Thus, we begin with some initial assumptions, ideas, fundamental beliefs, core values, axioms. We take certain things for granted. We try to make these as self-evident as possible.\footnote{  But of course, self-evidence, just like beauty, is in the eye of the beholder. There is more that can be said here.}  And from there we build our argument step by step, using logic as our guide. We define new constructs in terms of older, already accepted ones, and thus attempt to build a new world which has a solid foundation. 

Mathematical rationality can be found in the various versions of the ontological argument for the existence of God, as well as in the Declaration of Independence of the United States of America (see \cite{Grabiner} for a convincing argument about how mathematical rationality is built into the Declaration as well as many other illustrative examples of the impact of mathematical ways of knowing on Western thought). 

In the history of intellectual thought Rationalism of the European Enlightenment was met with a backlash movement, Romanticism. Today we can but do not have to see these two as directly opposing and mutually exclusive methodologies, each rejecting and invalidating the other. Alternatively, and I believe more productively, we might choose to accept that they point to two distinct ways of knowing, and occasionally certain truths will be more accessible via one way than another. 

Mathematical rationalism also has a similar complement, in what I will call mathematical imagination, or mathematical romanticism, if you will. Mathematical imagination is the way we select our axioms, the way we fix our principles. Mathematical imagination is how we determine our target truths. Human mathematicians do not start with a random formal axiomatic system and automatically go through all possible provable truths of the theory determined by it. Instead they engage with the worlds around them, both real and imaginary, and detect what is interesting, conjecture what might be productive to pursue, and then set out. Our human mathematics has a freedom to it and mathematical imagination captures this freedom. 

And freedom, broadly construed, may be viewed in these terms as well. S\'{a}ndor Szathm\'{a}ri's utopian, satiric novel, {\it The Voyage to Kazohinia}, might just be the best (fictional) guide to how mathematics is fundamental to a human society. Susan Siggelakis \cite{Siggelakis} describes how the protagonist of the novel learns that in a society without mathematics, ``nothing stable exists with which a human can connect and find meaning in his/her life''. There is only chaos and violence. As Edward Frenkel says in \cite{Frenkel} ``where there is no mathematics, there is no freedom''. 

\vspace{-10pt} \section*{mathematical ways of knowing: universals and eclecticisms}

Two other mathematical ways of knowing beckon us here: Universals and eclecticisms. The tendency of the mathematician to generalize is well known. ``Mathematics compares the most diverse phenomena and discovers the secret analogies that unite them,'' wrote Jean Baptiste Joseph Fourier. ``The art of doing mathematics is forgetting about the superfluous information,'' says Hendrik Lenstra. Thus the human mathematician tries to generalize, to abstract from specific examples, to reach universal statements with ``any'' and ``all'' that capture the essence of what is true about a whole slew of eclectic examples. The human desire to ``see the big picture'', the human tendency to ``find patterns'' is precisely what I mean by mathematical universalism. 

Accompanying and complementing (again productively) this tendency is the alternative, what I will call mathematical eclecticism. David Hilbert describes the complementarity of these two tendencies as follows:
\begin{quote}
In mathematics, as in any scientific research, we find two tendencies present. On the one hand, the tendency toward abstraction seeks to crystallize the logical relations inherent in the maze of material that is being studied, and to correlate the material in a systematic and orderly manner. On the other hand, the tendency toward intuitive understanding fosters a more immediate grasp of the objects one studies, a live rapport with them, so to speak, which stresses the concrete meaning of their relations.
\end{quote} 
Thus by mathematical eclecticism I mean the search for that one representative example on the one hand (``the art of doing mathematics consists in finding that special case which contains all the germs of generality,'' wrote Hilbert), and the excitement of the weirdness of eclectic cases on the other. In fact a lot of mathematics concerns itself with concrete examples. Paul Halmos wrote, ``the heart of mathematics consists of concrete examples and concrete problems.'' John B. Conway wrote, ``mathematics is a collection of examples; a theorem is a statement about a collection of examples and the purpose of proving theorems is to classify and explain the examples. . .'' 

``We think in generalities, but we live in details,'' said Alfred North Whitehead. In fact, I believe we think in both and we live in both. Once again, these two mathematical ways of knowing complement one another and help us live our human lives.

\vspace{-10pt} \section*{mathematical ways of knowing: certainty and ambiguity}

Alan Lightman says it best \cite{Lightman}: ``We are idealists and we are realists. We are dreamers and we are builders. We are experiencers and we are experimenters. We long for certainties, yet we ourselves are full of the ambiguities of the Mona Lisa and the I Ching. We ourselves are a part of the yin-yang of the world." Certainty and ambiguity find their way into mathematics and mathematical ways of knowing in a similarly complementary fashion. 

Certainty is a part of most people's perception of mathematics. Mathematics, for those people, is made up of questions and answers. Answers are certain once we know them. Indeed mathematics is perhaps the only certain knowledge we will ever have. I do not wish to minimize this perspective. I admit that I too have a romantic attachment to the idea that mathematical knowledge has a quality of certainty that goes farther than any other type of knowledge. And it took us a long time to get over this perspective as the unique way to conceive of mathematical truths and mathematical knowledge. 

The loss of certainty in mathematics began more than a century ago \cite{Kline}.  Lakatos was also influential in convincing many who cared to listen of the fallibility of mathematics \cite{Lak}. And Byers \cite{Byers} is perhaps the most detailed expositor of the role of ambiguity in the work of mathematics today. 

The complementarity of these two ways of knowing is rich and, at least to me, inspiring.

\vspace{-10pt} \section*{applications of mathematical ways of knowing}

Doing math at school or anywhere else is tied deeply in with our views of ourselves. This has good and bad aspects of course. We can relate our mathematical experiences to confidence, resilience, and determination as well as feelings of inadequacy, resistance, and rebelliousness; and as many can personally attest, any combination of the six can occur together. There is much emotion in mathematical engagement: hatred, love, anger, fear, anxiety, surprise, frustration, anticipation. How we handle mathematical challenges (suffering alone, valiantly standing defeated or undefeated, finding commiserators and conspirators) tells us about ourselves. Many of those who continue to do math after school connect with it at an emotional and personal level. We find aesthetic stimulation and creative joy in mathematical activity, as well as terrible frustration and occasional bouts of tedium. 

But can mathematical ways of knowing allow us to reach self-knowledge and a sense of identity mathematically? Andres Sanchez recounts how through an intentional application of set theory to his own personal life, he was able to discover his true identity and sexuality \cite{Sanchez}. Set theory, more generally, offers us pathways of thinking about belonging and not belonging. A clever student of mine, when asked to form study groups, chose to name his team ``the identity element'' as he wanted to work alone for the project in question. Mathematical ways of knowing may come in handy when thinking in terms of borders and boundaries of nations, communities, and cultural groups. They can help us think of inequalities and intersectionalities in different ways \cite{Cheng}. They can help us with moral dilemmas \cite{Voss, Voss2}. Mathematical ways of thinking and knowing can indeed allow us to view and understand ourselves in new and insightful ways.\footnote{  And in turn other ways of knowing can help us understand our mathematics better \cite{Gutierrez}.}

\vspace{-10pt} \section*{parting words -- till next time\ldots}

We have come a long way since Ptolemy argued that ``mathematics alone yields knowledge and that, furthermore, it is the only path to the good life'' \cite{Feke}. Mathematical knowledge has lost its certainty somewhere along the way, and primality in the eye of the public a long time ago. Mathematicians eventually learned to be humbler about the reaches of mathematical ways of knowing. But mathematics can still yield powerful knowledge, and not just the kind that can blow up cities and optimize factory production. I urge us to try and open up to the world once again. If we dig deeper into mathematical ways of knowing and the contexts where they might apply, mathematics may yet surprise us.


{
\setlength{\baselineskip}{13pt}

}


\begin{thebibliography}{99}


\bibitem{Anderson}
Anderson, Lorin W.; Krathwohl, David R., eds. (2001). A taxonomy for learning, teaching, and assessing: A revision of Bloom's taxonomy of educational objectives. Allyn and Bacon. 


\bibitem{Bloom}
Bloom, B. S.; Engelhart, M. D.; Furst, E. J.; Hill, W. H.; Krathwohl, D. R. (1956). Taxonomy of educational objectives: The classification of educational goals. Handbook I: Cognitive domain. New York: David McKay Company.

\bibitem{Byers} Byers, W. (2007). {\it How Mathematicians Think:\ using ambiguity, contradiction, and paradox to create mathematics}. Princeton:\ Princeton University Press.


\bibitem{Cheng} Eugenia Cheng, TEDx talk on how abstract mathematics can help us understand the world, available at \url{https://www.ted.com/talks/dr_eugenia_cheng_how_abstract_mathematics_can_help_us_understand_the_world}, 
last accessed on April 15, 2019. 

\bibitem{Feke}
Feke, Jacqueline (2018). Ptolemy's Philosophy: Mathematics as a Way of Life. Princeton University Press.


\bibitem{Frenkel}
Frenkel, Edward (2013). Love and Math: The Heart of Hidden Reality. Basic Books. 

\bibitem{Gettier} Edmund Gettier, ``Is Justified True Belief Knowledge?", {\it Analysis}, Volume {\bf 23} Issue 6 (1963), pages 121--123.


\bibitem{Grabiner}
Grabiner, Judith V. (1988). "The Centrality of Mathematics in the History of Western Thought," Mathematics Magazine, Volume 61 Issue 4, 220-230.

\bibitem{Grosholz}
Grosholz, Emily (2013).  ``Teaching the Complex Numbers: What History and Philosophy of Mathematics Suggest," Journal of Humanistic Mathematics, Volume 3 Issue 1, 62-73.

\bibitem{Harel}
Harel, Guershon (2008). ``What is Mathematics? A Pedagogical Answer to a Philosophical Question''. In Bonnie Gold and Roger Simons (eds.), Proof and Other Dilemmas: Mathematics and Philosophy. Mathematical Association of America. 265-290.

\bibitem{Gutierrez}
Guti\'{e}rrez, Rochelle (2012). ``Embracing Nepantla: Rethinking `Knowledge' and its Use in Mathematics Teaching'', REDIMAT - Journal of Research in Mathematics Education, Volume 1 Issue 1, 29-56. 

\bibitem{Hersh} Hersh, Reuben (2014). Experiencing Mathematics: What Do We Do, When We Do Mathematics? Providence: American Mathematical Society.

\bibitem{Karaali15}
Karaali, Gizem (2015). ``Can Zombies Do Math?'' In Mariana Bockarova, Marcel Danesi, Dragana Martinovic and Rafael Nunez (eds.), Mind in Mathematics: Essays on Mathematical Cognition and Mathematical Method, Interdisciplinary Studies on the Nature of Mathematics \#3, Lincom. 140?153. 

\bibitem{Kleiner}
Kleiner, Israel (1991). ``Rigor and Proof in Mathematics: A Historical Perspective'', Mathematics Magazine, Volume 64 Issue 5, 291-314.

\bibitem{Kline}
Kline, Morris (1980). Mathematics: The Loss of Certainty. Oxford: Oxford University Press. 

\bibitem{Lak}  Lakatos, I. (1976). {\it Proofs and Refutations:\ the logic of mathematical discovery}. New York:\ Cambridge University Press.

\bibitem{Lakoff} Lakoff, G. and N\'{u}\~{n}ez,  R. E. (2000). {\it Where Mathematics Comes From:\ How the Embodied Mind Brings Mathematics into Being}.  New York:\ Basic Books.


\bibitem{Lightman}
Lightman, Alan (2018). Searching for Stars on an Island in Maine. Penguin Random House. 

\bibitem{Pollard}
Pollard, Stephen (2014). ``Review of Experiencing Mathematics by Reuben Hersh'', Philosophia Mathematica, Volume 22 Issue 2, 271-274.

\bibitem{Sanchez}
Sanchez, Andres (2018). "My Sets and Sexuality," Journal of Humanistic Mathematics, Volume 8 Issue 1, 359-370. https://scholarship.claremont.edu/jhm/vol8/iss1/18/

\bibitem{Siggelakis}
Siggelakis, Susan J. (2019). ``Sándor Szathmári's Kazohinia: Mathematics and the Platonic Idea'', Journal of Humanistic Mathematics, Volume 9 Issue 1, 3-23. https://scholarship.claremont.edu/jhm/vol9/iss1/3/

\bibitem{Smith}
Smith, M.S., \& Stein, M.K. (1998). ``Selecting and creating mathematical tasks: From research to practice." Mathematics teaching in the middle school, Volume 3 Issue 5, 344-350

\bibitem{Tao} Terence Tao, ``There?s more to mathematics than rigour and proofs", blog pot available at \url{https://terrytao.wordpress.com/career-advice/theres-more-to-mathematics-than-rigour-and-proofs/}, last accessed on April 15, 2019.  

\bibitem{Tarski3356} 
Tarski, Alfred (1933/1956). ``The concept of truth in the languages of the deductive sciences'' (Polish), Prace Towarzystwa Naukowego Warszawskiego, Wydzial III Nauk Matematyczno-Fizycznych 34, Warsaw; reprinted in Zygmunt 1995, 13-172; expanded English translation in Tarski 1983 [1956], 152-278.

\bibitem{Tarski46} 
Tarski, Alfred (1946). Introduction to Logic and the Methodology of the Deductive Sciences. Oxford University Press.

\bibitem{Tarski56}
Tarski, A. and Vaught, R., 1956, ``Arithmetical extensions of relational systems'', Compositio Mathematica, Volume 13, 81-102.

\bibitem{Voss}
Sarah Voss, ``A Workshop to Introduce Concepts of Moral Math", {\it Journal of Humanistic Mathematics}, Volume {\bf 2} Issue 2 (July 2012), pages 114--128.  Available at: \url{https://scholarship.claremont.edu/jhm/vol2/iss2/10}, last accessed on April 15, 2019. 

\bibitem{Voss2}
Sarah Voss, ``Fuzzy Logic in Health Care Settings: Moral Math for Value-Laden Choices", {\it Journal of Humanistic Mathematics}, Volume {\bf 6} Issue 2 (July 2016), pages 161--178. Available at: \url{https://scholarship.claremont.edu/jhm/vol6/iss2/12}, last accessed on April 15, 2019. 














\end{thebibliography}
\end{document}